\documentclass[12pt]{amsart} 
\usepackage{amssymb, fullpage, color}

\newcommand{\be}{\begin{enumerate}}
\newcommand{\ee}{\end{enumerate}}

\DeclareMathOperator{\rank}{rank}

\renewcommand{\setminus}{-}  
\renewcommand{\epsilon}{\varepsilon}  
\newcommand{\vecalpha}{\vec{\alpha}}
\newcommand{\vecbeta}{\vec{\beta}}
\newcommand{\norm}{{\mathbf N}}
\newcommand{\bad}{{\operatorname{bad}}}

\newcommand{\Aff}{{\mathbb A}}

\newcommand{\F}{{\mathbb F}}

\newcommand{\N}{{\mathbb N}}
\newcommand{\Q}{{\mathbb Q}}
\newcommand{\R}{{\mathbb R}}
\newcommand{\Z}{{\mathbb Z}}
\newcommand{\PP}{{\mathbb P}}
\newcommand{\Qbar}{{\overline{\Q}}}

\newcommand{\kbar}{{\overline{k}}}
\newcommand{\Kbar}{{\overline{K}}}

\newcommand{\Fbar}{{\overline{\F}}}

\newcommand{\pp}{{\mathfrak p}}
\newcommand{\qq}{{\mathfrak q}}

\newcommand{\calE}{{\mathcal E}}

\newcommand{\calL}{{\mathcal L}}
\newcommand{\calM}{{\mathcal M}}

\newcommand{\calP}{{\mathcal P}}

\newcommand{\calS}{{\mathcal S}}
\newcommand{\calT}{{\mathcal T}}
\newcommand{\calU}{{\mathcal U}}

\newcommand{\calW}{{\mathcal W}}

\newcommand{\calZ}{{\mathcal Z}}

\newcommand{\OO}{{\mathcal O}}

\DeclareMathOperator{\End}{End}

\DeclareMathOperator{\Lie}{Lie}

\DeclareMathOperator{\Aut}{Aut}
\DeclareMathOperator{\Gal}{Gal}

\DeclareMathOperator{\ord}{ord}

\DeclareMathOperator{\Frob}{Frob}
\DeclareMathOperator{\disc}{disc}

\newcommand{\tors}{{\operatorname{tors}}}

\newcommand{\GL}{\operatorname{GL}}

\newcommand{\isom}{\simeq}

\newcommand{\Union}{\bigcup}
\newcommand{\union}{\cup}

\newcommand{\directsum}{\oplus}

\newtheorem{theorem}{Theorem}[section]
\newtheorem{lemma}[theorem]{Lemma}
\newtheorem{corollary}[theorem]{Corollary}
\newtheorem{proposition}[theorem]{Proposition}

\theoremstyle{definition}
\newtheorem{definition}[theorem]{Definition}
\newtheorem{question}[theorem]{Question}

\theoremstyle{remark}
\newtheorem{remark}[theorem]{Remark}

\begin{document}

\title[Diophantine definability]{Diophantine definability of infinite
discrete non-archimedean sets and Diophantine models over large subrings
of number fields}
\subjclass[2000]{Primary 11U05; Secondary 11G05}
\keywords{Hilbert's Tenth Problem, elliptic curve, Mazur's Conjecture, Diophantine definition, norm equation}

\author{Bjorn Poonen}
\address{Department of Mathematics, University of California,
Berkeley, CA 94720-3840, USA}
\email{poonen@math.berkeley.edu}
\urladdr{http://math.berkeley.edu/\~{}poonen}

\author{Alexandra Shlapentokh} \address{Department of Mathematics, East Carolina University, Greenville, NC 27858,
USA} \email{shlapentokha@mail.ecu.edu} \urladdr{http://www.personal.ecu.edu/shlapentokha} \thanks{B.P. was
partially supported by a Packard Fellowship and by NSF grants DMS-9801104 and DMS-0301280. A.S. was partially
supported by NSF grants DMS-9988620 and DMS-0354907, and by ECU Faculty Senate Summer 2003 Grant.}

\date{August 19, 2004}

\begin{abstract}
We prove that infinite $\pp$-adically discrete sets have Diophantine
definitions in large subrings of some number fields. First, if $K$ is
a totally real number field or a totally complex degree-$2$ extension
of a totally real number field, then there exists a prime $\pp$ of $K$
and a set of $K$-primes $\calS$ of density arbitrarily close to $1$
such that there is an infinite $\pp$-adically discrete set that is
Diophantine over the ring $\OO_{K,\calS}$ of $\calS$-integers in
$K$. Second, if $K$ is a number field over which there exists an
elliptic curve of rank $1$, then there exists a set of $K$-primes
$\calS$ of density $1$ and an infinite Diophantine subset of
$\OO_{K,\calS}$ that is $v$-adically discrete for every place $v$ of
$K$. Third, if $K$ is a number field over which there exists an
elliptic curve of rank $1$, then there exists a set of $K$-primes
$\calS$ of density $1$ such that there exists a Diophantine model of
$\Z$ over $\OO_{K,\calS}$. This line of research is motivated by a
question of Mazur concerning the distribution of rational points on
varieties in a non-archimedean topology and questions concerning
extensions of Hilbert's Tenth Problem to subrings of number fields.
\end{abstract}

\maketitle

\section{Introduction}
\setcounter{equation}{0}

Matijasevi\v{c} (following work of Davis, Putnam, and Robinson)
proved that Hilbert's Tenth Problem could not be solved:
that is, there does not exist an algorithm, that given
an arbitrary multivariable polynomial equation $f(x_1,\dots,x_n)=0$ with
coefficients in $\Z$, decides whether or not a solution in $\Z^n$ exists.
It is not known whether an analogous algorithm exists, however,
if in the problem one replaces $\Z$ by $\Q$ in both places.
One natural approach to proving a negative answer for $\Q$
is to show that $\Z$ admits a Diophantine definition over $\Q$,
or more generally that there is a Diophantine model of the ring $\Z$
over $\Q$;
the meaning of these statements is given in
Definitions \ref{D:diophantine definition} and~\ref{D:diophantine model}
below.
(See~\cite{hilbertstenthproblem} for an introduction to the subject.)

\begin{definition}
\label{D:diophantine definition}
Let $R$ be a (commutative) ring.
Suppose $A \subseteq R^k$ for some $k \in \N$.
Then we say that $A$ has a {\em Diophantine definition over $R$}
if there exists a polynomial
\[
    f(t_1,\ldots,t_k, x_1,\ldots,x_n) \in R[t_1,\ldots,t_k,,x_1,\ldots,x_n]
\]
such that for any $(t_1,\ldots,t_k) \in R^k$,
\[
    (t_1,\ldots,t_k) \in A
        \quad\iff\quad
    \exists x_1,\ldots,x_n \in R, \; f(t_1,\ldots,t_k,x_1,...,x_n) = 0.
\]
In this case we also say
that $A$ is a {\em Diophantine subset} of $R^k$,
or that $A$ is {\em Diophantine over~$R$}.
\end{definition}

\begin{remark}
Suppose that $R$ is a domain whose quotient field
is not algebraically closed.
Then
\begin{enumerate}
\item[(a)]
Relaxing Definition~\ref{D:diophantine definition}
to allow an arbitrary finite conjunction of equations
in place of the single equation on the right hand side
does not enlarge the collection of Diophantine sets.
\item[(b)]
Unions and intersections of Diophantine sets
are Diophantine.
\end{enumerate}
See the introduction in \cite{pheidas1994} for details.
\end{remark}

\begin{definition}
\label{D:diophantine model}
A {\em Diophantine model of $\Z$ over a ring $R$} is a Diophantine
subset $A \subseteq R^k$ for some $k$
together with a bijection $\phi\colon \Z \to A$
such that the graphs of addition and multiplication
(subsets of $\Z^3$) correspond under $\phi$ to Diophantine subsets
of $A^3 \subseteq R^{3k}$.
\end{definition}

Mazur formulated a conjecture that would imply that a
Diophantine definition of $\Z$ over $\Q$ does not exist,
and later in ~\cite{cornelissen-zahidi2000}
it was found that his conjecture also
ruled out the existence of a Diophantine model of $\Z$ over $\Q$.
One form of Mazur's conjecture was that for a variety $X$ over $\Q$,
the closure of $X(\Q)$ in the topological space $X(\R)$ should have
at most finitely many connected components.
See \cite{mazur1992}, \cite{mazur1994},  \cite{mazur1995},
\cite{mazur1998}, \cite{ct-skorobogatov-sd1997},
\cite{cornelissen-zahidi2000}, \cite{poonensubrings},
and \cite{shlapentokh2003}
for more about the conjecture and its consequences.

Mazur also formulated an analogue applying to both archimedean and
nonarchimedean completions of arbitrary number fields.
Specifically, on page 257 of~\cite{mazur1998} he asked:

\begin{question}
\label{ques:1}
Let $V$ be any variety defined over a number field $K$.
Let $\calS$ be a finite set of places of $K$,
and consider $K_{\calS}= \prod_{v \in \calS} K_v$ viewed as
locally compact topological ring.
Let $V(K_{\calS})$ denote the topological space
of $K_{\calS}$-rational points.
For every point $p \in V(K_{\calS})$
define $W(p) \subset V$ to be the subvariety defined over $K$
that is the intersection of Zariski closures of the subsets $V(K) \cap U$,
where $U$ ranges through all open neighborhoods of $p$ in $V(K_{\calS})$.
As $p$ ranges through the points of $V(K_{\calS})$,
are there only a finite number of distinct subvarieties $W(p)$?
\end{question}

In Question~\ref{ques:1}, it does not matter whether we require
$V$ to be irreducible.
(We will not.)

\begin{proposition}
Fix a number field $K$ and a place $\pp$.
If Question~\ref{ques:1} has a positive answer for $K$ and $\calS:=\{\pp\}$,
then there does not exist
an infinite, $\pp$-adically discrete, Diophantine subset of $K$.
\end{proposition}

\begin{proof}
Suppose there exists a subset $A$ of $K$ that is
infinite, $\pp$-adically discrete, and Diophantine over $K$.
The Diophantine definition of $A$ corresponds to an affine algebraic set $V$
such that the projection $\pi: V \to \Aff^1$ onto the first coordinate
satisfies $\pi(V(K))=A$.

Suppose $a \in A$.
Since $A$ is $\pp$-adically discrete,
there exists an open neighborhood $N$ of $a$ in $K_\pp$
such that $A \cap N = \{a\}$.
Pick $p \in V(K)$ with $\pi(p)=a$.
Then $U:=\pi^{-1}(N)$ is an open neighborhood of $p$,
and $\pi$ maps $V(K) \cap U$ into $A \cap N = \{a\}$,
so $W(p) \subseteq \pi^{-1}(a)$.

By choosing one $p$ above each $a \in A$,
we get infinitely many disjoint subvarieties $W(p)$,
contradicting the positive answer to Question~\ref{ques:1}.
\end{proof}

\bigskip

In view of the proposition above,
constructing a Diophantine definition of
an infinite discrete $\pp$-adic set over a number field $K$ would be one way
to answer Question \ref{ques:1} (negatively) for $K$.
Unfortunately, at the moment such a construction seems out of reach.
Thus instead we consider analogues in which $K$
is replaced by some of its large integrally closed subrings.

\begin{definition}
If $K$ is a number field, let $\calP_K$ be the set of finite primes of $K$.
For $\calS \subseteq \calP_K$, define the {\em ring of $\calS$-integers}
\[
    \OO_{K,\calS} =
    \{\, x \in K \mid
    \ord_\pp x \geq 0 \text{ for all $\pp \notin {\calS}$}\,\}.
\]
(Elsewhere the term $\calS$-integers often presupposes that $\calS$ is finite,
but we will use this term for infinite $\calS$ also.)
\end{definition}

If $\calS = \emptyset$,
then $\OO_{K,\calS}$ equals the ring $\OO_K$ of algebraic integers of $K$.
If $\calS = \calP_K$, then $\OO_{K,\calS}=K$.
In general, $\OO_{K,\calS}$ lies somewhere between $\OO_K$ and $K$,
and the density of $\calS$ (if it exists)
may be used as a measure of the ``size'' of $\OO_{K,\calS}$.
Throughout this paper, ``density'' means {\em natural density},
which is defined as follows.

\begin{definition}
\label{def:natural}
Let $\calS \subseteq \calP_K$.
The {\em density} of $\calS$ is defined to be the limit
\[
    \lim_{X\rightarrow \infty}
    \frac{\#\{\pp \in \calS: N\pp \leq X\}}
    {\#\{\mbox{all } \pp: N\pp \leq X\}}
\]
if it exists.
If the density of $\calS$ does not exist,
we can replace the limit in Definition \ref{def:natural}
by $\limsup$ or $\liminf$ and hence define
the upper or lower densities of $\calS$.
(See \cite[VIII,~\S4]{lang1994} for more about density.)

\end{definition}

The study of Diophantine definability and the archimedean conjecture
of Mazur over rings of $\calS$-integers
has produced Diophantine definitions of $\Z$
and discrete archimedean sets over large subrings of some number fields
(see  \cite{shlapentokh1997}, \cite{shlapentokh2000highdensity},
\cite{shlapentokh2002}, \cite{shlapentokh2004} and \cite{shlapentokh2003}).
Recently in \cite{poonensubrings},
the first author constructed an infinite discrete
Diophantine set (in the archimedean topology) and a
Diophantine model of $\Z$ over a subring of $\Q$
corresponding to a set of primes of density~$1$.
Thus he showed that the analogue of Hilbert's Tenth Problem
is undecidable over such a ring.

In this paper we consider Diophantine definability of
infinite discrete $\pp$-adic sets over some rings of $\calS$-integers.
Our results will come from two sources:
norm equations (as in \cite{shlapentokh2003})
and elliptic curves (as in \cite{poonensubrings}).
Our main results are stated below.
When we say that a subset of $\OO_{K,\calS}$ is Diophantine,
we mean that it is Diophantine {\em over $\OO_{K,\calS}$}.
A subset $\calS$ of $\calP_K$ is {\em recursive}
if there exists an algorithm that takes as input an element of $K$
(given by its coordinates with respect to some fixed $\Q$-basis)
and decides whether it belongs to $\OO_{K,\calS}$.
 
\begin{theorem}
\label{thm:mainnorm}
Let $K$ be a totally real number field or
a totally complex degree-$2$ extension of a totally real number field.
Let $\pp$ be any prime of $K$. Then for any $\varepsilon >0$
there exists a recursive
set of $K$-primes $\calS$ containing
$\pp$ of density $> 1-\varepsilon$ such that there exists an infinite Diophantine subset of $\OO_{K,\calS}$ that
is discrete and closed when viewed as a subset of the completion $K_\pp$.
\end{theorem}

\begin{theorem}
\label{thm:main2}
Let $K$ be any number field for which there exists
an elliptic curve $E$ such that $\rank E(K)=1$.
Then
\begin{enumerate}
\item
There exist recursive subsets $\calT_1,\calT_2 \subseteq \calP_K$
of density~$0$ such that for any $\calS$ with
$\calT_1 \subseteq \calS \subseteq \calP_K-\calT_2$,
there exists an infinite Diophantine subset $A$ of $\OO_{K,\calS}$
such that for all places $v$ of $K$,
the set $A$ is discrete when viewed as a subset of the completion $K_v$.
\item
There exist recursive subsets $\calT_1',\calT_2' \subseteq \calP_K$
of density~$0$ such that for any $\calS$ with
$\calT_1' \subseteq \calS \subseteq \calP_K-\calT_2'$,
there exists a Diophantine model of the ring $\Z$ over $\OO_{K,\calS}$.
\end{enumerate}
\end{theorem}

\section{Using Norm Equations}
\setcounter{equation}{0}

In this section we use norm equations to
construct infinite Diophantine $\pp$-adically discrete sets,
in order to prove Theorem~\ref{thm:mainnorm}.

\subsection{Preliminary results}

\begin{proposition}
\label{prop:notzero}
Let $K$ be a number field and let $\calS \subseteq \calP_K$.
The $\OO_{K,\calS} - \{0\}$ is Diophantine over $\OO_{K,\calS}$.
\end{proposition}

\begin{proof}
See Proposition~2.6 on page 113 of~\cite{shlapentokh2000survey}.
(Note: there is a typo in the statement in~\cite{shlapentokh2000survey}:
the last $K$ should be $\OO_{K,W}$.)
\end{proof}

The importance of Proposition~\ref{prop:notzero} is that
equations with variables intended to range over $K$
can now be interpreted in the arithmetic of $\OO_{K,\calS}$,
since elements of $K$ can be represented as fractions
of elements of $\OO_{K,\calS}$ with nonzero denominator.

\begin{proposition}
\label{prop:dvr}
Let $K$ be a number field, and let $\pp \in \calP_K$.
Then the discrete valuation ring
$\OO_{K,\calP_K - \{\pp\}} = \{\, x\in K : \ord_\pp x \ge 0\,\}$
is Diophantine over $K$.
\end{proposition}

\begin{proof}
See Lemma~3.22 in~\cite{shlapentokh1994}.
\end{proof}

\begin{corollary}
\label{cor:order relations}
Let $K$ be a number field, and let $\pp \in \calP_K$.
Then the sets
$\{\, x \in K : \ord_\pp x > 0\,\}$,
$\{\,(x,y) \in K^2 : \ord_\pp x \ge \ord_\pp y \,\}$,
and
$\{\,(x,y) \in K^2 : \ord_\pp x = \ord_\pp y \,\}$
are Diophantine over $K$.
\end{corollary}

\begin{proof}
Fix $a \in K$ with $\ord_\pp a=1$.
Then
\begin{align*}
    \ord_\pp x > 0 \quad&\iff\quad \ord_\pp(x/a) \ge 0, \\
    \ord_\pp x \ge \ord_\pp y
    \quad&\iff\quad
    (\exists r) (x=ry \;\;\;\text{and}\;\; \ord_\pp r \ge 0), \\
    \ord_\pp x = \ord_\pp y
    \quad&\iff\quad
    (\ord_\pp x \ge \ord_\pp y)  \text{ and } (\ord_\pp y \ge \ord_\pp x).
\end{align*}
\end{proof}

Propositions \ref{prop:notzero} and~\ref{prop:dvr}
together imply the following generalization of Proposition~\ref{prop:dvr}
(cf.~Theorem 4.4 in~\cite{shlapentokh1994}):

\begin{proposition}
\label{prop:dedekind}
Let $K$ be a number field.
If $\calS \subseteq \calS' \subseteq \calP_K$
and $\calS' - \calS$ is finite,
then $\OO_{K,\calS}$ is Diophantine over $\OO_{K,\calS'}$.
\end{proposition}

\begin{lemma}
\label{le:max}
Let $F$ be a number field.
Let $\{\omega_1,\dots,\omega_s\}$ be a $\Z$-basis for $\OO_F$.
Let $a_1,\dots,a_s \in \Q$, and let $x=\sum_{i=1}^s a_i \omega_i$.
Let $p$ be a prime of $\Q$ that does not ramify in $F$.
Then
\[
    \min_i \ord_{p} a_i = \min_{\mathfrak P} \ord_{\mathfrak P} x,
\]
where $\mathfrak P$ ranges over $F$-primes above $p$.
\end{lemma}

\begin{proof}
Since $p$ is unramified in $F$, the ideal $p \OO_F$ factors as the
product of the $\mathfrak P$. Thus we have an equality
$p^m \OO_F = \prod_{\mathfrak P} {\mathfrak P}^m$ for any $m \in \Z$.
The two minimums in the statement equal the largest $m$ for which $x$ belongs
to this fractional ideal,
since $p^m \omega_1$, \dots, $p^m \omega_s$ form a $\Z$-basis for $p^m \OO_F$.
\end{proof}

\begin{lemma}
\label{L:E exists}
For any rational primes $p$ and $q$,
there exists a degree-$p$ cyclic extension $E/\Q$
in which $q$ splits completely.
If moreover $p$ is odd, then any such $E$ is totally real.
\end{lemma}

\begin{proof}
Choose a rational prime $\ell$ splitting completely
in $\Q(e^{2\pi i/p}, q^{1/p})$.
Then $\ell \equiv 1 \pmod{p}$ and the image of $q$
in $\left(\Z/\ell\Z\right)^* \isom \Gal(\Q(e^{2\pi i/\ell})/\Q)$
is a $p$-th power.
Equivalently, the Frobenius automorphism
$\Frob_q \in \Gal(\Q(e^{2\pi i/\ell})/\Q)$
belongs to the index-$p$ subgroup $H$ of
the cyclic group $\Gal(\Q(e^{2\pi i/\ell})/\Q)$.
Let $E$ be the subfield of $\Gal(\Q(e^{2\pi i/\ell})/\Q)$ fixed by $H$.
Then the image of $\Frob_q$ in $\Gal(E/\Q)$ is trivial,
so $q$ splits completely in $E$.
Odd-degree Galois extensions of $\Q$ are totally real.
\end{proof}

\subsection{Notation and Assumptions}
\label{subsec:assumptions}

We view all number fields as being subfields
of a fixed algebraic closure $\Qbar$,
so that compositums are well-defined.

\begin{itemize}
\item Let $K$ be the number field given in Theorem~\ref{thm:mainnorm}. Thus $K$ is totally real, or $K$ is a
totally complex degree-$2$ extension of a totally real field. \item Let $n=[K:\Q]$. \item Let $\pp$ be the prime
of $K$ in Theorem~\ref{thm:mainnorm}. \item Let $p_{\Q}$ be the prime of $\Q$ below $\pp$. \item Let
$p=p_0<p_1<\cdots<p_n$ be a sequence of odd primes such that $p_i>n$ and $1/p_i < \epsilon/(n+1)$ for all $i$.
\item Let $E=E_0$, $E_1$, \dots, $E_n$ be a sequence of totally real cyclic extensions of $\Q$ such
that $[E_i:\Q]=p_i$ and $p_\Q$ splits completely in $E$. (These exist by Lemma~\ref{L:E exists}.) \item Let $L$ be
an imaginary degree-$2$ extension of $\Q$ in which $p_{\Q}$ splits completely. (For instance, let $b<0$ be an
integer that is a nonzero square mod $p_{\Q}$ and let $L=\Q(\sqrt{b})$.) \item Let $\mathcal V_{\Q}$ be the set of
rational primes that are inert in all of the extensions $E_i/\Q$ for $0 \le i \le n$. Let ${\mathcal W}_{\Q} =
{\mathcal V}_{\Q} \cup \{p_{\Q}\}$.
\item Let ${\mathcal W}_{EL}$ be the set of $EL$-primes above ${\mathcal W}_{\Q}$.
\item Let ${\mathcal W}_{K}$ be the set of $K$-primes above ${\mathcal W}_{\Q}$. \item Let $\sigma_L$, $\sigma_E$
be generators of $\Gal(L/\Q)$, $\Gal(E/\Q)$ respectively. \item Let $\pp_{EL}$ be a prime above $p_{\Q}$ in
$EL$.
\item Let $\Omega=\{\omega_1,\ldots,\omega_{2p}\}$ be a $\Z$-basis for $\OO_{EL}$. Then $\Omega$ is also an
$\OO_{\Q,\calW_\Q}$-basis for $\OO_{EL,\calW_{EL}}$.
\end{itemize}

\subsection{Discrete Diophantine subsets of large subrings of $\Q$}

The following proposition provides the foundation for our construction.

\begin{proposition}
\label{prop:normsolutions}
Let $x \in \OO_{EL, \calW_{EL}}$ be a solution to the system
\begin{equation}
\label{sys:norm}
\left \{ \begin{array}{c} N_{EL/L}(x)=1,\\
    N_{EL/E}(x) = 1.\end{array} \right .
\end{equation}
such that
\begin{gather}
\label{it:1}
    \text{$x$ is a unit at all primes of $EL$ above $p_{\Q}$
    except possibly for } \\
\nonumber
    \text{
    $\pp_{EL}$,
    $\sigma_E(\pp_{EL})$,
    $\sigma_L(\pp_{EL})$, and
    $\sigma_E\sigma_L(\pp_{EL}) = \sigma_L\sigma_E(\pp_{EL})$; and}
\end{gather}
\begin{equation}
\label{it:2}
    \text{$\ord_{\pp_{EL}}x >0$
    and $\ord_{\sigma_E\sigma_L(\pp_{EL})}x >0$.}
\end{equation}
Then the divisor of $x$ equals
\begin{equation}
\label{eq:divisor}
\left(\frac{\pp_{EL}\,\sigma_E\sigma_L(\pp_{EL})}
{\sigma_E(\pp_{EL})\,\sigma_L(\pp_{EL})}\right)^k
\end{equation}
for some $k \in \Z_{>0}$.
Furthermore, the set of solutions
to \eqref{sys:norm}
satisfying \eqref{it:1} and~\eqref{it:2}
is nonempty and closed under multiplication.
If $x_1$ and $x_2$ are solutions to \eqref{sys:norm},
then $x_1/x_2$ is again a solution to \eqref{sys:norm},
but $x_1/x_2$ will not necessarily satisfy \eqref{it:1} and~\eqref{it:2}
even if $x_1$ and $x_2$ do.
\end{proposition}

\begin{proof}
Since $p_\Q$ splits completely in $E$ and in $L$,
it splits completely in $EL$.
Now apply Section~3.3.2 on page~131 of \cite{shlapentokh2000survey}.
\end{proof}
 
\begin{corollary}
\label{cor:Dioph}
Let $A$ be the set of
$(a_1, \ldots, a_{2p})$ in $\left(\OO_{\Q,\calW_\Q}\right)^{2p}$
such that
the element $x=\sum_{i=1}^{2p} a_i \omega_i$ of $\OO_{EL,\calW_{EL}}$
satisfies \eqref{sys:norm}, \eqref{it:1}, and~\eqref{it:2}.
Let $B$ be the set of $b \in \OO_{\Q,\calW_\Q}$
such that for some $(a_1,\dots,a_{2p}) \in A$,
\begin{enumerate}
\item[(i)] The element $b$ equals one of the $a_i$, and
\item[(ii)]
    $\ord_{p_{\Q}} b =
    \min\{\ord_{p_{\Q}}a_1,\ldots,\ord_{p_{\Q}}a_{2p}\}.$
\end{enumerate}
(Thus $B$ is the set of ``$p_\Q$-adically largest coordinates''
of elements of $A$.)
For $r \in \Z$, let $B_r = \{\,b \in B : \ord_{p_\Q} b = r \,\}$.
Then
\be
\item $A$ and $B$ are Diophantine over $\OO_{\Q,\calW_\Q}$.
\item
There exists $M \in \Z_{>0}$ such that
$B = \Union_{m=1}^\infty B_{-Mm}$
and each $B_{-Mm}$ in the union is a nonempty finite set.
\ee
\end{corollary}

\begin{proof}\hfill
\be
\item
By writing each norm in~\eqref{sys:norm} as a product of conjugates,
we find that the equations~\eqref{sys:norm} are equivalent
to polynomial equations in the $a_i$.
By Corollary~\ref{cor:order relations}
together with the ``Going up and then down'' method
(see Section 2.2 of \cite{shlapentokh2000survey}),
conditions \eqref{it:1} and~\eqref{it:2} are Diophantine.
Thus $A$ is Diophantine.

Condition~(i) is $\prod_{i=1}^{2p}(b-a_i)=0$,
and condition~(ii) is Diophantine by Corollary~\ref{cor:order relations},
so $B$ is Diophantine too.
\item
If $b \in B$ is associated to $(a_1,\dots,a_{2p}) \in A$
and $k$ is the positive integer in~\eqref{eq:divisor}
for the element $x = \sum_{i=1}^{2p} a_i \omega_i$,
then $\ord_{p_\Q} b = -k$ by Lemma~\ref{le:max}.
Thus the set of $r$ for which $B_r \ne \emptyset$
equals the set of possibilities for $-k$ in~\eqref{eq:divisor}
as $x$ varies over all solutions to
\eqref{sys:norm} satisfying \eqref{it:1} and~(\ref{it:2}).
Let $M$ be the smallest positive integer such that $-M$
is a possible value of $-k$.
Proposition~\ref{prop:normsolutions}
implies that the set of possibilities for $-k$
is then $\{-M,-2M,\dots\}$.
Thus $B_r \ne \emptyset$ if and only if $r \in \{-M,-2M,\dots\}$.

It remains to prove that each $B_r$ is finite.
Fix $r \in \Z$.
The finiteness of $B_r$ follows
once we show that the set of $x$
satisfying \eqref{sys:norm}, \eqref{it:1}, and~(\ref{it:2})
and having $k=-r$ in~\eqref{eq:divisor}
is finite.
Suppose $x_1$ and $x_2$ are two such values of $x$,
so by~\eqref{eq:divisor} the divisor of $x_1/x_2$ is trivial.
Then $x_1/x_2 \in \OO_{EL}^*$ and $N_{EL/E}(x_1/x_2)=1$.
But $EL/E$ is a totally complex degree-$2$ extension
of a totally real field, so these conditions imply
that $x_1/x_2$ is a root of unity.
Finally, there are only finitely many roots of unity in $EL$.
\ee
\end{proof}

\begin{corollary}
\label{cor:discrete for Q}
There exists an infinite Diophantine subset $B$ of $\OO_{\Q,\calW_\Q}$
that is $p_{\Q}$-adically discrete and closed.
\end{corollary}

\begin{proof}
The Diophantine set $B$ of Corollary~\ref{cor:Dioph} is
infinite, because it is an infinite disjoint union of nonempty sets.
It is $p_{\Q}$-adically discrete and closed,
because for each $r<0$,
the set of $b \in B$ with $\ord_{p_\Q} b \ge r$ is finite.
\end{proof}

\subsection{Discrete Diophantine subsets of large subrings of $K$}

\begin{proposition}
\label{prop:A}
Let $\calU$ be a set of $K$-primes remaining prime in $E_iK/K$
for $i=0,\dots,n$.
Then
there exists a set of $K$-primes $\bar{\calU}$ such that
the set difference
$({\calU} \setminus \bar{\calU} )\cup
(\bar{\calU} \setminus {\calU})$
is finite and $\OO_{K,\bar{\calU}} \cap \Q$ has a Diophantine definition
over $\OO_{K,\bar{\calU}}$.
\end{proposition}

\begin{proof}
See Corollary~2.3 and Theorem~3.8 of~\cite{shlapentokh2002}.
(The application of Corollary~2.3 requires that $p_i$ does not
divide the absolute degree of the Galois closure of $K/\Q$;
this holds since $p_i>n$.)
\end{proof}

Clearly Proposition~\ref{prop:A} implies the slightly stronger version
in which the hypothesis on $\calU$ is weakened to the
hypothesis that {\em all but finitely many}
primes of $\calU$ are inert in all the $E_iK/K$.
We now show that we can also insist that the new set of primes contains the
original set:

\begin{proposition}
\label{prop:B}
Let $\calU$ be a set of $K$-primes such that all but finitely many
of them are inert in $E_iK/K$ for $i=0,\dots,n$.
Then
there exists a set of $K$-primes ${\calU}'$ {\em containing $\calU$} such that
${\calU}' - \calU$ is finite
and
$\OO_{K,{\calU}'} \cap \Q$ is Diophantine over $\OO_{K,{\calU}'}$.
\end{proposition}

\begin{proof}
Let $\bar{\calU}$ be the set given by
(the slightly stronger version of) Proposition~\ref{prop:A},
and let ${\calU}' = \calU \cup \bar{\calU}$.
Thus ${\calU}' - \calU$ is finite.

By choice of $\bar{\calU}$,
the set $R:=\OO_{K,\bar{\calU}} \cap \Q$ is Diophantine over
$\OO_{K,\bar{\calU}}$,
which is Diophantine over $\OO_{K,{\calU}'}$
by Proposition~\ref{prop:dedekind}.
Thus $R$ is a Diophantine subset of $\OO_{K,{\calU}'}$.
The desired subset $\OO_{K,{\calU}'} \cap \Q$
can now be defined as the set of elements of $\OO_{K,{\calU}'}$
equal to a ratio of elements of $R$ with nonzero denominator
(here we use Proposition~\ref{prop:notzero} for $R$).
\end{proof}

\begin{proof}[Proof of Theorem~\ref{thm:mainnorm}]
We will apply Proposition~\ref{prop:B} to the set $\calW_K$
defined in Section~\ref{subsec:assumptions}.
By assumption,
the cyclic extension $E_i/\Q$
has prime degree $p_i > n$.
Thus $[E_iK:K]=p_i$,
and moreover, a $K$-prime is inert in $E_iK/K$
if and only if the $\Q$-prime below it is inert in $E_i/\Q$.
All primes of $\calW_\Q$ but $p_\Q$ are inert in all the $E_i/\Q$,
so all primes of $\calW_K$ but finitely many are inert in all the $E_iK/K$.
Thus we may apply Proposition~\ref{prop:B} to find a set of $K$-primes
$\calW'_K$ such that $\calW'_K - \calW_K$ is finite
and such that the set $\OO_{\Q,\calW'_\Q}$ is Diophantine over
$\OO_{K,\calW'_K}$,
where $\calW'_\Q$ is the set of $\Q$-primes $q$ such that
all $K$-primes above $q$ lie in $\calW'_K$.

Since $\calW'_K - \calW_K$ is finite,
$\calW'_\Q - \calW_\Q$ is finite.

Now the infinite set $B$ of Corollary~\ref{cor:discrete for Q}
is Diophantine over $\OO_{\Q,\calW_\Q}$,
which by Proposition~\ref{prop:dedekind}
is Diophantine over $\OO_{\Q,\calW'_\Q}$,
which is Diophantine over $\OO_{K,\calW'_K}$.
Thus $B$ is Diophantine over $\OO_{K,\calW'_K}$.
Since $\pp$ lies above $p_\Q$,
the set $B$ is $\pp$-adically discrete and closed.
To fulfill the requirements of Theorem~\ref{thm:mainnorm},
we take $\calS = \calW'_K$.
 This is recursive, since up to a finite set, its primes
are characterized by splitting behavior in a finite list
of extension fields.

It remains to show that $\calW'_K$ has density greater than $1-\epsilon$.
Up to finitely many primes, $\calW'_K$
is defined by the splitting behavior in finitely many extensions of $K$
(namely, the $E_iK/K$).
Thus, by the Chebotarev Density Theorem
(see Th\'eor\`eme~1 of~\cite{serre1981} for a version using natural density),
$\calW'_K$ has a density.
The density of the set of $K$-primes that {\em fail} to be inert
in $E_iK/K$ is $1/p_i$,
so the density of $\calW'_K$
is at least
\[
    1 - \sum_{i=0}^n \frac{1}{p_i}
    \;>\; 1 - \sum_{i=0}^n \frac{\epsilon}{n+1}
    \;=\; 1-\epsilon.
\]
\end{proof}

\section{Using Elliptic Curves}
\setcounter{equation}{0}

\subsection{Notation}
\label{S:notation section}
\begin{itemize}
\item Whenever $k$ is a perfect field, let $\kbar$ be an algebraic closure,
and let $G_k=\Gal(\kbar/k)$.
\item $K$ is a number field.
\item $E$ is an elliptic curve of rank~1 over $K$.
(In particular, we assume that $K$ is such that such an $E$ exists).
\item We fix a Weierstrass equation $y^2=x^3+ax+b$ for $E$
with coefficients in the ring of integers of $K$.
\item $E(K)_\tors$ is the torsion subgroup of $E(K)$.
\item $r$ is an even multiple of $\#E(K)_\tors$.
\item $Q \in E(K)$ is such that $Q$ generates $E(K)/E(K)_\tors$.
\item $P:=rQ$.
\item $\calP_\Q=\{2,3,5,\dots\}$ is the set of rational primes.
\item $\calP_K$ is the set of all finite primes of $K$.
\item Let ${\mathcal S}_\bad \subseteq \calP_K$ consist of the primes
that ramify in $K/\Q$, the primes for which the reduction
of the chosen Weierstrass model is singular
(this includes all primes above $2$),
and the primes at which the coordinates of $P$ are not integral.
We occasionally view $E$ as the scheme over $\OO_{K,\calS_\bad}$
defined by the homogenization of the Weierstrass equation.
\item $E'$ is the smooth affine curve $y^2=x^3+ax+b$ over $\OO_{K,\calS_\bad}$.
Thus $E'$ is $E$ with the zero section removed.
\item ${\calM}_K$ is the set of all normalized absolute values of $K$.
\item ${\calM}_{K,\infty} \subset {\calM}_K$ is the set of all archimedean absolute values of $K$.
\item For $\pp \in \calP_K$, let
    \be
        \item $K_{\pp}$ be the completion of $K$ at $\pp$.
        \item $R_\pp$ be the valuation ring of $K_\pp$
        \item $\F_{\pp}$ be the residue field of $R_\pp$,
        \item $\norm\pp = \#\F_\pp$ be the absolute norm of $\pp$
    \ee

\item Write $nP=(x_n,y_n)$ where $x_n,y_n \in K$.
\item Let the divisor of $x_n$ be of the form
\[
    \frac{{\mathfrak a}_n}{{\mathfrak d}_n}{\mathfrak b}_n,
\]
where
\begin{itemize}
\item ${\mathfrak d}_n=\prod_{{\mathfrak q}}{\mathfrak q}^{-a_{{\mathfrak q}}}$, where the product is taken over
all primes ${\mathfrak q}$ of $K$ not in $\calS_{\bad}$ such that $a_{{\mathfrak q}}=\ord_{{\mathfrak q}}x_n <0$,
\item ${\mathfrak a}_n=\prod_{{\mathfrak q}}{\mathfrak q}^{a_{{\mathfrak q}}}$, where the product is taken over
all primes ${\mathfrak q}$ of $K$ not in $\calS_{\bad}$ such that $a_{{\mathfrak q}}=\ord_{{\mathfrak q}}x_n >0$.
\item ${\mathfrak b}_n = \prod_{{\mathfrak q}}{\mathfrak q}^{a_{{\mathfrak q}}}$, where the product is taken over
all primes ${\mathfrak q}\in \calS_{\bad}$ and  $a_{{\mathfrak q}}=\ord_{{\mathfrak q}}x_n$.
\end{itemize}
\item For $n$ as above, let ${\calS}_n = \{\pp \in {\mathcal P}_K : \pp | {\mathfrak d}_n\}$. By definition of
$\calS_\bad$, we have $\calS_1=\emptyset$.
\item For $n$ as above, let $d_n = \norm {\mathfrak d}_n \in \Z_{\ge 1}$.
\item For $u \in K^*$ and $v$ a place of $K$ lying above the place $p$ of $\Q$ (possibly $p=\infty$), define the
(unnormalized) local height $h_v(u) = \log \max\{ \|u\|_v,1\}$ where $\|u\|_v = |\norm_{K_v/\Q_p}(u)|_p$ and $K_v$
and $\Q_p$ denote completions.
\item For $u \in K^*$, define the global height $h(u)=\sum_{v \in {\mathcal M}_K} h_v(u)$.
\item For $\ell \in \calP_\Q$, define $a_\ell$ to be the smallest positive number such that ${\mathcal
S}_{\ell^{a_\ell}} \ne \emptyset$. (By Siegel's Theorem, $a_\ell=1$ for all but finitely many $\ell$.)
\item Let $\calL = \{\ell\in {\mathcal P}_\Q: a_\ell >1\}$
and $L = \prod_{\ell \in \calL} \ell^{a_\ell-1}$.
\item Let $\pp_\ell$ be a prime of largest norm
in $\calS_{\ell^{a_\ell}}$.
\item For $\ell,m \in \calP_\Q$,
let $\pp_{\ell m}$ be a prime of largest norm in
$\calS_{\ell m} \setminus (\calS_\ell \cup \calS_m)$,
if this set is nonempty (see Proposition~\ref{prop:biggerS}).
\item $E[m]$ denotes $\{\,T \in E(\Kbar) : mT=0\,\}$.
\item If $\lambda$ is an ideal in the endomorphism ring of $E$,
then
\[
    E[\lambda] := \{\, T \in E(\Kbar) :
    aQ=0 \text{ for some $a \in \lambda$} \,\}.
\]
 \end{itemize}

\subsection{Divisibility of denominators of $x$-coordinates}

The next lemma is the number field analogue to
Lemma~3.1(a) of \cite{poonensubrings},
where it was proved for $\Q$.

\begin{lemma}
\label{le:Po3.1}
Let $P \in E(K)\setminus\{0\}$ be of infinite order
and let  $n \in \Z\setminus\{0\}$.
Let ${\mathfrak r}$ be an integral divisor of $K$.
Then $\{n \in \Z: {\mathfrak r} \mid {\mathfrak d}_n(P)\}$
is a subgroup of $\Z$.
\end{lemma}

\begin{proof}
It is enough to prove the lemma when $\mathfrak r$
is a prime power $\pp^m$.
Let $\hat{E}$ be the formal group over $R_\pp$
defined by the chosen Weierstrass model of $E$.
Let $\hat{E}(\pp R_\pp)$ be the group of points associated to $\hat{E}$.
There exist Laurent series $x(z)=z^{-2}+\cdots$ and $y(z)=z^{-3}+\cdots$
with coefficients in $R_\pp$
giving an injective homomorphism $\hat{E}(\pp R_\pp) \to E(K_\pp)$
whose image is the set $E_1(K_\pp)$ of
$(x,y) \in E(K_\pp)$ with $\ord_\pp(x)<0$
(together with $O$):
this follows from Proposition~VII.2.2 of~\cite{silvermanAEC}
when the Weierstrass equation is minimal,
but the proof there does not use the minimality.
Since $\ord_\pp x(z) = -2 \ord_\pp z$ whenever $z \in \pp R_\pp$,
the set of $(x,y) \in E(K_\pp)$ with $\ord_\pp(x) \le -m$
(together with $O$)
corresponds under this homomorphism to
a subgroup $\pp^{\lceil m/2 \rceil} R_\pp$ of $\pp R_\pp$,
and hence is a subgroup of $E(K_\pp)$.
\end{proof}

\begin{corollary}
\label{cor:intersec}
Let $m,n \in \Z - \{0\}$,
and let $(m,n)$ be their gcd.
Then ${\mathcal S}_m \cap {\mathcal S}_n = {\mathcal S}_{(m,n)}$.
In particular, if $(m,n)=1$ then
${\mathcal S}_m \cap {\mathcal S}_n=\emptyset$.
\end{corollary}

\subsection{New primes in denominators of $x$-coordinates}

\begin{lemma}
\label{le:orderchange}
Let $n \in \Z_{\ge 1}$.
Suppose that ${\mathfrak t} \in \calP_K$ divides ${\mathfrak d}_n$,
and $p \ge 3$ is a rational prime.
\be
 \item If ${\mathfrak t} \mid p$,
then $\ord_{{\mathfrak t}}{\mathfrak d}_{pn}=
\ord_{{\mathfrak t}}{\mathfrak d}_n+2$.
 \item If ${\mathfrak t} \nmid p$,
then $\ord_{{\mathfrak t}}{\mathfrak d}_{pn}=
\ord_{{\mathfrak t}}{\mathfrak d}_n$.
\ee
\end{lemma}

\begin{proof}
 We will use the notation of Lemma \ref{le:Po3.1}.
Since ${\mathfrak t} \mid {\mathfrak d}_n$,
by assumption ${\mathfrak t} \not \in \calS_\bad$.
In particular $\mathfrak t$ is not ramified over $\Q$.
Furthermore, $\ord_{{\mathfrak t}} x_n <0$,
so $nP \in E_1(K_{{\mathfrak t}})$.
Let $z$ be the corresponding element in the group of points
$\hat{E}({\mathfrak t} R_{\mathfrak t})$ of the formal group.
We have $[p]z = pf(z)  + g(z^p)$, by Proposition 2.3,
page 116 and Corollary 4.4, page 120 of~\cite{silvermanAEC},
where  $[p]$ is the multiplication-by-$p$ in the formal group
and $f(T), g(T) \in R_{\mathfrak t}[[T]]$
satisfy $g(0) =0$ and $f(T)=T + \mbox{ higher order terms}$.
Thus $\ord_{\mathfrak t}([p]z)$ equals $(\ord_{\mathfrak t} z) + 1$
or $\ord_{\mathfrak t} z$,
depending on whether ${\mathfrak t} \mid p$.
Since $x=z^{-2} + \cdots$,
we find that
$\ord_{\mathfrak t} x_{pn}$ equals $(\ord_{\mathfrak t} x_n) - 2$
or $\ord_{\mathfrak t} x_n$,
depending on whether ${\mathfrak t} \mid p$.
\end{proof}

\begin{lemma}
 \label{le:denomheight}
There exists $c \in \R_{>0}$
such that $\log d_n =(c-o(1))n^2$ as $n \longrightarrow \infty$.
\end{lemma}

\begin{proof} Let $\hat{h}$ be the canonical height on $E(K)$.
Then $\hat{h}(nP)/n^2$ is a positive constant
independent of $n$. The Weil height differs from $\hat{h}$ by $O(1)$,
so $h(x_n)/n^2$ tends to a positive limit as $n \to \infty$.
By definition, $h(x_n)$ differs from $\log d_n$ by the
sum of $h_v(x_n)$ over archimedean $v$ and $v \in \calS_\bad$.
By the theorem on page 101 of~\cite{serremordellweil},
$h_v(x_n)/h(x_n) \to 0$ as $n \to \infty$ for each $v$,
so $(\log d_n)/h(x_n)$ tends to $1$ as $n \to \infty$.
Thus $(\log d_n)/n^2$ tends to a positive limit as $n \to \infty$.
\end{proof}

The next proposition is a number field version of
Lemma~3.4 of \cite{poonensubrings}.

\begin{proposition}
\label{prop:biggerS}
 If $\ell, m \in \calP_\Q$ and $\max(\ell,m)$ is sufficiently large,
then $\calS_{\ell m} \setminus (\calS_\ell \cup \calS_m) \ne \emptyset$.
\end{proposition}

\begin{proof}
Suppose $\calS_{\ell m} \setminus (\calS_\ell \cup \calS_m) = \emptyset$.
We claim that
${\mathfrak d}_{\ell m} \mid \ell^2 m^2 {\mathfrak d}_\ell {\mathfrak d}_m$.
To check this, we compare orders of both sides at a prime
$\mathfrak t$ dividing ${\mathfrak d}_{\ell m}$.
By assumption, $\mathfrak t$ divides either ${\mathfrak d}_\ell$ or
${\mathfrak d}_m$.
Without loss of generality, assume ${\mathfrak t} \mid {\mathfrak d}_\ell$.
Then Lemma~\ref{le:orderchange}
implies
$\ord_{\mathfrak t} {\mathfrak d}_{\ell m}
\le \ord_{\mathfrak t}(m^2 {\mathfrak d}_{\ell})$,
which proves the claim.

Taking norms, we obtain $d_{\ell m} \mid (\ell m)^{2[K:\Q]} d_\ell d_m$.
Taking logs and applying Lemma~\ref{le:denomheight},
we deduce
\begin{align*}
    (c-o(1)) \ell^2 m^2
    &\le 2[K:\Q] (\log \ell + \log m) + (c-o(1)) \ell^2 +(c-o(1)) m^2 \\
    & \le (c+o(1))(\ell^2 + m^2),
\end{align*}
which is a contradiction once $\ell$ and $m$ are both sufficiently large.
\end{proof}

\subsection{Density of prime multiples}

\begin{lemma}
\label{L:dirichlet-vinogradov}
Let $\vecalpha \in \R^n$,
let $I$ be an open neighborhood of $0$ in $\R^n/\Z^n$,
and let $d \in \Z_{\ge 1}$.
Then the set of primes $\ell \equiv 1 \pmod{d}$
such that $(\ell-1)\vecalpha \mod 1$ is in $I$
has positive upper density.
\end{lemma}

\begin{proof}
Let $\vecalpha=(\alpha_1,\dots,\alpha_n)$.
We first reduce to the case that $1,\alpha_1,\dots,\alpha_n$
are $\Z$-independent.
Choose $m \in \Z_{\ge 1}$ and $\beta_1,\dots,\beta_r$
such that $1/m,\beta_1,\dots,\beta_r$ form a $\Z$-basis
for the subgroup of $\R$ generated by $1,\alpha_1,\dots,\alpha_n$.
Replacing $d$ by a positive integer multiple only reduces the density,
so we may assume $m \mid d$.
For fixed $i$, if $\alpha_i = \frac{c_0}{m} + \sum_{j=1}^r c_j \beta_j$,
then for $\ell \equiv 1 \pmod{d}$ we have
\[
    (\ell-1) \alpha_i \equiv \sum_{j=1}^r c_j (\ell-1)\beta_j \pmod{1},
\]
so it suffices to prove positivity of the upper density of
$\ell \equiv 1 \pmod{d}$ for which $(\ell-1)\beta_j$
is sufficiently close mod~$1$ to $0$ for all $j$.

In fact, we will prove the stronger result that
the points $(\ell-1)\vecbeta \in (\R/\Z)^r$
for prime $\ell \equiv 1 \pmod{d}$
are equidistributed.
By Weyl's equidistribution criterion \cite[Satz~3]{weyl1916},
we reduce to proving that for any $\alpha \in \R \setminus \Q$,
\[
    \sum_{\substack{\ell \le x \\ \ell \equiv 1 \pmod{d}}}
    e^{2\pi i \ell \alpha}
    = o(\pi(x))
\]
as $x \to \infty$.
This is a consequence of Vinogradov's work on exponential sums over primes:
see \cite[p.~34]{montgomery1994}, for instance.
\end{proof}

\subsection{Denominators of $x$-coordinates having many small prime factors}

We next prove an analogue of Lemma~7.1 of~\cite{poonensubrings}
showing that it is rare that $\calS_\ell$ has a large fraction
of the small primes.
For prime $\ell$, define
\[
    \mu_\ell = \sup_{X \in \Z_{\ge 2}}
        \frac{\# \{\pp \in \calS_\ell : \norm\pp \le X\}}
        {\# \{\pp \in \calP_K: \norm\pp \le X\}}.
\]

\begin{lemma}
\label{L:mulemma} For any $\epsilon>0$, the density of $\{\, \ell : \mu_\ell > \epsilon \,\}$ is $0$.
\end{lemma}

\begin{proof}
The proof can be copied from that of Lemma~7.1 of~\cite{poonensubrings},
using ``the primes $\pp \in \calP_K$ with $\norm\pp \le X$''
everywhere in place of the ``the primes $p$ up to $X$'':
it requires only the facts
\begin{enumerate}
\item
The function $\pi_K(X):=\# \{\pp \in \calP_K : \norm\pp \le X\}$
is $(1+o(1))X/\log X$ as $X \to \infty$ (Theorem~3 on page 213
of~\cite{cassels-frohlich}).
\item
$\#\calS_\ell \le \log_2 d_\ell$
(clear, since each prime has norm at least $2$)
\item
$\log_2 d_\ell = O(\ell^2)$ as $\ell \to \infty$
(follows from Lemma~\ref{le:denomheight}).
\end{enumerate}
\end{proof}

\subsection{Construction of the $\ell_i$}
\label{S:construction}

By \cite[Corollary~VI.5.1.1]{silvermanAEC}
and \cite[Corollary~V.2.3.1]{silvermanATAEC},
there is an isomorphism of real Lie groups
$\prod_{v \in \calM_{K,\infty}} E(K_v) \isom (\R/\Z)^N \times (\Z/2\Z)^{N'}$
for some $N \ge 1$ and $N' \ge 0$.
Fix such an isomorphism, and embed $E(K)$
diagonally in $\prod_{v \in \calM_{K,\infty}} E(K_v)$.
Since $P=rQ$ with $r$ even,
the point $P$ maps to an element $\vecalpha \in (\R/\Z)^{N}$.

Define a sequence of primes $\ell_i$ inductively as follows.
Given $\ell_1,\dots,\ell_{i-1}$,
let $\ell_i$ be the smallest prime outside $\calL$
and exceeding the bound implicit in Proposition~\ref{prop:biggerS}
such that all of the following hold:
\begin{enumerate}
\item $\ell_i>\ell_j$ for all $j<i$,
\item $\mu_{\ell_i} \le 2^{-i}$,
\item $\norm\pp_{\ell_i \ell_j} > 2^i$ for all $j<i$,
\item $\norm\pp_{\ell \ell_i} > 2^i$ for all $\ell \in \calL$,
\item $\ell_i \equiv 1 \pmod{i!}$, and
\item $|x_{\ell_i-1}|_v > i$ for all $v \in \calM_{K,\infty}$.
\end{enumerate}

\begin{proposition}
\label{P:well-defined}
The sequence $\ell_1,\ell_2,\dots$ is well-defined and computable.
\end{proposition}

\begin{proof}
Condition~(6) is equivalent to the requirement that
$(\ell-1)\vecalpha$ lie in a certain open neighborhood of $0$ in
$(\R/\Z)^N$, since the Lie group isomorphism maps neighborhoods of $O$
to neighborhoods of $0$.  Thus by Lemma~\ref{L:dirichlet-vinogradov},
the set of primes satisfying (5) and~(6) has positive upper
density. By Lemma~\ref{L:mulemma}, (2) fails for a set of density
$0$. Therefore it will suffice to show that (1), (3), and~(4) are
satisfied by all sufficiently large $\ell_i$.

For fixed $j \le i$, the primes $\pp_{\ell_i \ell_j}$ for varying
values of $\ell_i$ are distinct by Corollary~\ref{cor:intersec},
so eventually their norms are greater than $2^i$.
The same holds for $\pp_{\ell \ell_i}$ for fixed $\ell \in \calL$.
Thus by taking $\ell_i$ sufficiently large,
we can make all the $\pp_{\ell_i \ell_j}$ and $\pp_{\ell \ell_i}$
have norm greater than $2^i$.
Thus the sequence is well-defined.

Each $\ell_i$ can be computed by searching primes in increasing order
until one is found satisfying the conditions:
condition~(6) can be tested effectively,
since $|x_{\ell_i-1}|_v$ is an algebraic real number.
\end{proof}

As in Section~4 of~\cite{poonensubrings},
we define the following subsets of $\calP_K$:
\begin{itemize}
\item $\calT_1 = \calS_\bad \union \Union_{i \ge 1} \calS_{\ell_i}$,
\item $\calT_2^a$ is the set of $\pp_\ell$
    for $\ell \notin \{\ell_1,\ell_2,\dots\}$,
\item $\calT_2^b = \{\, \pp_{\ell_i \ell_j} : 1 \le j \le i \,\}$,
\item $\calT_2^c = \{\, \pp_{\ell \ell_i} : \ell \in \calL, i \ge 1 \,\}$, and
\item $\calT_2 = \calT_2^a \union \calT_2^b \union \calT_2^c$.
\end{itemize}

As in Section~5 of~\cite{poonensubrings},
we prove that

\begin{lemma}
\label{L:integer points}
The sets $\calT_1$ and $\calT_2$ are disjoint.
If the subset $\calS \subset \calP_K$
contains $\calT_1$ and is disjoint from $\calT_2$,
then $\calE := E'(\OO_{K,\calS}) \cap rE(K)$ is the union of
$\{\, \pm \ell_i P: i \ge 1 \,\}$
and some subset of the finite set
$\left\{\, sP : s \mid \prod_{\ell \in \calL} \ell^{a_\ell-1} \,\right\}$.
\end{lemma}

\begin{proof}
Once we note that $rE(K)=\Z P$,
the proofs proceed as in Section~5 of~\cite{poonensubrings}.
\end{proof}

The recursiveness of $\calT_1$ and $\calT_2$
follows as in Section~8 of~\cite{poonensubrings}, using the following:

\begin{lemma}
\label{exactorder}
If $\ell$ is prime, then $\ell \mid \#E(\F_{\pp_\ell})$.
\end{lemma}

\begin{proof}
For $\pp \notin \calS_\bad$,
a multiple $nP$ reduces to $0$ in $E(\F_{\pp})$ if and only if
$\pp$ divides ${\mathfrak d}_n$.
Hence, by definition of $\pp_\ell$,
the point $\ell^{a_\ell} P$ reduces to $0$ in
$E(\F_{\pp_\ell})$ but $\ell^{a_\ell-1} P$ does not.
\end{proof}

\subsection{Density of $\calT_1$ and $\calT_2$}

The proofs that $\calT_1$, $\calT_2^b$, and $\calT_2^c$ have density $0$
are identical to the proofs in Section~9 of~\cite{poonensubrings}.
The remainder of this section is devoted to proving that
$\calT_2^a$ has density $0$.
Again, we follow~\cite{poonensubrings}, but more work is necessary because
we no longer assume that $E$ has no CM.

For $n \in \Z_{>0}$,
let $\omega(n)$ be the number of distinct prime factors of $n$.

\begin{lemma}
For any $t \ge 1$,
the density of $\{\,\pp : \omega(\#E(\F_\pp))<t \,\}$ is $0$.
\end{lemma}

\begin{proof}
If $E$ does not have CM (i.e., $\End E_\Kbar = \Z$), then the proof
given for $K=\Q$ in Lemma~9.3 of~\cite{poonensubrings} generalizes
easily to arbitrary $K$. The first step in this proof, which we will
also use for the CM case, is to relate divisibility to Galois
representations: namely, for a prime $\pp$ of good reduction not above
$\ell$, the condition $\ell \mid \#E(\F_\pp)$ is equivalent
to the existence of an $\ell$-torsion point in $E(\Fbar_\pp)$
fixed by the Frobenius element of $\Gal(\Fbar_\pp/\F_\pp)$,
which in turn is equivalent to the condition that the
the image of a Frobenius element at $\pp$ under
$G_K \to \Aut E[\ell] = \GL_2(\Z/\ell\Z)$
has a nontrivial fixed vector in $(\Z/\ell\Z)^2$.

We assume from now on that $E$ has CM,
say by an order $\OO$ in a quadratic imaginary field $F$.
Using the action
of $\OO$ on $\Lie E$ (the tangent space of $E$ at $O$),
we may view $\OO$ (and hence also $F$) as a
subring of $\Kbar$.

All the endomorphisms are defined over
the compositum $KF$ \cite[II.2.2(b)]{silvermanATAEC}.

Let $\Lambda$ be the set of primes $\ell$ of $\Z$
such that $\ell \nmid \disc(\OO)$ and $(\ell)$ factors into
distinct prime ideals $\lambda$ and $\bar{\lambda}$ of $\OO$.
(Later we will delete finitely many primes from $\Lambda$.)
 
For $\ell \in \Lambda$,
\[
     E[\ell] = E[\lambda] \directsum E[\bar{\lambda}].
\]

The summands on the right are
free modules over $\OO/\lambda$ and $\OO/\bar{\lambda}$,
respectively,
and we choose generators for each in order to obtain identifications
\[
E[\ell] \isom \OO/\lambda \directsum \OO/\bar{\lambda} \isom (\Z/\ell\Z)^2,
\]
 and hence $\Aut E[\ell] \isom \GL_2(\Z/\ell\Z)$.
The action of $G_{KF}$ commutes with the $\OO$-action,
so the image of $G_{KF}$ in
$\GL_2(\Z/\ell\Z)$ lies in the subgroup
$\begin{pmatrix} * & \\ & * \end{pmatrix}$ of diagonal matrices.
On the other hand,
if $\tau \in G_K - G_{KF}$,
then $\tau$ interchanges $\lambda$ and $\bar{\lambda}$,
and its image in $\GL_2(\Z/\ell\Z)$
lies in the coset $\begin{pmatrix} & * \\ * & \end{pmatrix}$.

Define subgroups
\begin{align*}
    D &= \begin{pmatrix} * &  \\  & * \end{pmatrix} \\
    H &= \begin{pmatrix} * &  \\  & * \end{pmatrix}
        \cup \begin{pmatrix}  & * \\ * &  \end{pmatrix}
\end{align*}
of $\GL_2 \left(\prod_{\ell \in \Lambda} \Z/\ell\Z \right)$.
Thus the image of
$\rho \colon G_K \to \GL_2 \left(\prod_{\ell \in \Lambda} \Z/\ell\Z \right)$
lies in $H$.
It is a classical fact that $\rho(G_{KF})$ is open in $D$
(see for example, the Corollaire on page~302 of~\cite{serre1972}),
so $\rho(G_K)$ is open in $H$.
By deleting a few primes from $\Lambda$,
we may assume that $\rho(G_K)$ contains $D$,
and hence equals $D$ or $H$,
depending on whether $F \subseteq K$ or not.

Let $\pi_\ell$ be the probability that a random element of the subgroup
$\begin{pmatrix} * & \\ & * \end{pmatrix} \subseteq \GL_2(\Z/\ell\Z)$
has a nontrivial fixed vector.
A calculation shows that $\pi_\ell = 2/\ell + O(1/\ell^2)$.
Since $\Lambda$ has density $1/2$, the series
$\sum_{\ell \in \Lambda} 1/\ell$ diverges.
Thus $\sum_{\ell \in \Lambda} \pi_\ell$ diverges.
Elementary probability shows that if $X_1,X_2,\dots$
are independent events, and the sum of their probabilities diverges,
then as $C \to \infty$,
the probability that fewer than $t$ of the first $C$ events occur
tends to $0$.
Therefore as $C \to \infty$,
if $\sigma$ is chosen uniformly at random from the image of $D$
in $\prod_{\ell \in \Lambda, \ell<C} \GL_2(\Z/\ell\Z)$,
then the probability that $\sigma$ has fewer than $t$ components
with a nontrivial fixed vector tends to $0$.

Similarly, a random element of the coset
$\begin{pmatrix}  & * \\ * &  \end{pmatrix} \subseteq \GL_2(\Z/\ell\Z)$
has a nontrivial fixed vector with probability
$1/\ell + O(1/\ell^2)$.
Thus the probability that a random $\sigma$ from the image
of  $H \setminus D$ in $\prod_{\ell \in \Lambda, \ell<C} \GL_2(\Z/\ell\Z)$
has fewer than $t$ components with a nontrivial fixed vector
tends to $0$.

Combining the previous two paragraphs shows that the same holds
for a random element of the image $I_C$ of
$G_K \to \prod_{\ell \in \Lambda, \ell<C} \GL_2(\Z/\ell\Z)$.
But the Chebotarev Density Theorem
implies that the images of the Frobenius elements $\Frob_\pp$
in $I_C$ are equidistributed in $I_C$.
Finally we apply the ``key first step'' mentioned
at the beginning of this proof to get the desired result.
\end{proof}

Using the preceding lemma,
the proof that $\calT_2^a$ has
follows exactly the proof of Proposition~9.4 in~\cite{poonensubrings},
using $\norm\pp$ in place of $p$ in the inequalities.

Thus $\calT_1$ and $\calT_2$ have density $0$.

\subsection{Convergence and discreteness}

\begin{lemma}
\label{L:convergence}
For each $v \in \calM_K$,
the sequence $\ell_1 P$, $\ell_2 P$, $\dots$ converges in $E(K_v)$ to $P$.
\end{lemma}

\begin{proof} It suffices to show that $(\ell_i-1)P \to O$
in $E(K_v)$ as $i \to \infty$. If $v$ is
archimedean, this holds by condition~(6) in the construction of the $\ell_i$.
Now suppose $v$ is nonarchimedean.
Then the topological group $E(K_v)$ has a basis consisting
of open finite-index subgroups $U$,
namely the groups in the filtration
appearing in the proof of \cite[VII.6.3]{silvermanAEC}.
So it suffices to show, given $U$,
that $(\ell_i-1)P \in U$ for sufficiently large $i$.
Let $j$ be the index of $U$ in $E(K_v)$.
If $i \ge j$, then $j \mid i! \mid \ell_i-1$,
by condition~(5) in the construction of the $\ell_i$,
so $(\ell_i-1)P \in U$.

\end{proof}

\begin{proposition}
\label{P:Diophantine-discrete}
Let $\calS$ be as in Lemma~\ref{L:integer points}.
Let $A:=\{x_{\ell_1},x_{\ell_2},\dots\}$.
Then $A$ is a Diophantine subset of $\OO_{K,\calS}$.
For any $v \in \calM_K$,
the set $A$ is discrete when viewed as a subset of $K_v$.
\end{proposition}

\begin{proof}
By Lemma~\ref{L:integer points}, $x(\calE)$ is the union of
the set $A:=\{x_{\ell_1},x_{\ell_2},\dots\}$ and a finite set.
Since $\calE$ is Diophantine over $\OO_{K,\calS}$,
so is $A$.

By Lemma~\ref{L:convergence},
the elements of $A$ form a convergent sequence in $K_v$,
and the limit $x_1$ of the sequence is not in $A$,
so $A$ is discrete.
\end{proof}

This completes the proof of part~(1) of Theorem~\ref{thm:main2}.

\subsection{A Diophantine model of $\Z$}

We next show how to find a Diophantine model of the ring $\Z$
over certain rings $\OO_{K,\calS}$.

\begin{lemma}
Let $B=\{\,2^n+n^2 : n \in \Z_{\ge 1}\,\}$.
Multiplication admits a positive existential definition
in the structure $\calZ:=(\Z_{\ge 1},1,+,B)$.
(Here $B$ is considered as a unary predicate.)
\end{lemma}

\begin{proof}
We can define $>$ by
\[
    x > y \quad\iff\quad (\exists z) \; x=y+z
\]
and for fixed $a \in \Z$, we have
\[
    x \ne a \quad\iff\quad (x>a) \lor (a>x),
\]
so this predicate is positive existential in $\calZ$.
For fixed $c \in \Z_{\ge 1}$, the function $x \mapsto cx$ is
positive existential, since it can be obtained by repeated addition.

Call $x,y$ {\em consecutive}
if there exists $n \in \Z_{\ge 1}$ such that
$x=2^n+n^2$ and $y=2^{n+1}+(n+1)^2$.
The set of such $(x,y)$ is positive existential in $\calZ$
since it differs from
\[
    \{\,(x,y) \in B^2 : x < y < 3x\,\}
\]
in a finite set.
Next
\[
    \{\,((2y-z)-(2x-y),2x-y):
    \text{ $x,y$ are consecutive and $y,z$ are consecutive}\,\}
\]
equals the set $T:=\{\,(2n-1,n^2-2n-1) : n \in \Z_{\ge 1}\,\}$.
We have
\[
    (u=v^2) \wedge (v>0) \quad\iff\quad (2v-1,u-2v-1) \in T.
\]
Call this relation $P(u,v)$.
Then
\begin{align*}
    u=v^2 \quad&\iff\quad P(u,v) \vee P(u,-v) \vee ((u=0) \wedge (v=0)), \\
    u=vw \quad&\iff\quad (v+w)^2 = v^2 + w^2 + 2u,
\end{align*}
so we can construct a positive existential definition of multiplication.
\end{proof}

\begin{remark}
Y.~Matijasevi{\v{c}} (private communication)
independently discovered a recursive set $B$
such that multiplication admits a positive existential definition
in $(\Z_{\ge 1},1,+,B)$.

\end{remark}

\begin{corollary}
\label{C:positive existential model}
The structure $(\Z,0,1,+,\cdot)$ admits a positive existential model
in the structure $\calZ$.
\end{corollary}

Because of Corollary~\ref{C:positive existential model},
instead of finding a Diophantine model of the ring $\Z$ over $\OO_{K,\calS}$,
it will suffice to find a Diophantine model of $\calZ$.

Now we redo the construction in Section~\ref{S:construction},
but change some of the conditions defining the
sequence of primes $\ell_i$.
Fix $\pp,\qq \in \calP_K -\calS_\bad$ of degree~1
such that neither $\pp$ nor $\qq$ divides $y_1=y(P)$,
and such that the underlying primes $p,q \in \calP_\Q$ are distinct and odd.
Let $M=p q \#E(\F_\pp) \#E(\F_\qq)$.
Keep conditions (1) through~(4), but replace conditions (5) and~(6) by the
following:
\begin{itemize}
\item[($5'$)] $\ell_i \equiv 1 \pmod M$,
\item[($6'$)] the highest power of $p$
dividing $(\ell_i-1)/M$ is $p^i$, and
\item[($7'$)] $q$ divides $(\ell_i-1)/M$ if and only if $i \in B$.
\end{itemize}

\begin{proposition}
\label{P:well-defined2}
The sequence $\ell_1,\ell_2,\dots$ is well-defined and computable.
\end{proposition}

\begin{proof}
The proof is the same as that of Proposition~\ref{P:well-defined}
except that instead of Lemma~\ref{L:dirichlet-vinogradov},
we will use just Dirichlet's theorem on primes in arithmetic progressions
to show that conditions ($5'$) through~($7'$) are satisfied by a
positive density of primes $\ell_i$.

The conditions amount to congruences modulo $p^{i+1} q M$,
so it suffices to show that one of the congruence classes
is of the form $a \pmod{p^{i+1} q M}$ with $(a,p^{i+1} q M)=1$.
Define
\[
    a = \begin{cases}
        1+p^i q M & \text{if $i \in B$} \\
        1+(p^i q + p^{i+1}) M & \text{if $i \notin B$.}
    \end{cases}
\]
Then $a$ is congruent to $1$ modulo $p$, modulo $q$, and modulo $M$,
so $(a,p^{i+1} q M)=1$.
Any prime in the residue class $a \pmod{p^{i+1} q M}$
satisfies ($5'$) through~($7'$).
\end{proof}

\begin{lemma}
\label{L:xdifference}
If $m \in \Z_{\ge 1}$, then
\[
    \ord_\pp(x_{mM+1}-x_1) = \ord_\pp(x_{M+1}-x_1) + \ord_p m.
\]
\end{lemma}

\begin{proof} Let $R$ be the valuation ring $R_\pp$
defined in Section~\ref{S:notation section}.
Because $y_1 \in R^*$ and $\pp \nmid 2$,
for any $r \ge 1$, the restriction of
the $x$-coordinate map $E(R/\pp^r) \to \PP^1(R/\pp^r)$
to the subset of points of $E(R/\pp^r)$
with the same image in $E(R/\pp)$ as $P$
is injective (the $y$-coordinate can be recovered from the
the $x$-coordinate as the square root of an element of $R^*$:
its sign is determined by the fact that the point is in the residue
class of $P$).
Thus for $r \ge 1$, the ideal $\pp^r$ divides
$x_{mM+1}-x_1$ if and only if $(mM+1)P$ and $P$
have the same image in $E(R/\pp^r)$, or equivalently if $(mM)P$
maps to $O$ in $E(R/\pp^r)$.
Let $z \in \hat{E}(\pp R)$ be the point of the formal group
corresponding to $MP$.
Then the condition that $(mM)P$ maps to $O$ in $E(R/\pp^r)$
is equivalent to $[m](z) \in \pp^r$, where $[m]$
denotes the multiplication-by-$m$ map in the formal group.

It remains to prove that $\ord_\pp [m](z) = \ord_\pp z + \ord_p m$.
By induction on $m$, it suffices to prove this when $m$ is prime.
Then, in the proof of Lemma~\ref{le:orderchange},
we have $[m](z) = m f(z) + g(z^m)$
where $f(T),g(T) \in R[[T]]$ satisfy $g(0)=0$
and $f(T)=T+\text{higher order terms}$.
Finally $\ord_\pp m = \ord_p m$
is $1$ or $0$ according to whether $m=p$ or not,
so the result follows.
\end{proof}

\begin{proposition}
\label{P:model of calZ}
Let $\calS$ be as in Lemma~\ref{L:integer points}.
Let $A:=\{x_{\ell_1},x_{\ell_2},\dots\}$.
Then $A$ is a Diophantine model of $\calZ$ over $\OO_{K,\calS}$,
via the bijection $\phi\colon \Z_{\ge 1} \to A$ taking $i$ to $x_{\ell_i}$.
\end{proposition}

\begin{proof}
The set $A$ is Diophantine over $\OO_{K,\calS}$
by the argument in the proof of Proposition~\ref{P:Diophantine-discrete}.

We have
\begin{align*}
    i \in B
    \quad&\iff\quad \text{$q$ divides $(\ell_i-1)/M$}
            &&\text{(by condition~($7'$))}\\
    \quad&\iff\quad \ord_\qq(x_{\ell_i}-x_1) > \ord_\qq(x_{M+1}-x_1),
\end{align*}
by Lemma~\ref{L:xdifference} (with $\qq$ in place of $\pp$).
The latter inequality is a Diophantine condition on $x_{\ell_i}$,
by Corollary~\ref{cor:order relations}
(in which we represent elements of $K$ as ratios of
elements of $\OO_{K,\calS}$).
Thus the subset $\phi(B)$ of $A$ is Diophantine over $\OO_{K,\calS}$.

Finally, for $i \in \Z_{\ge 1}$,
Lemma~\ref{L:xdifference} and condition~($6'$)
imply $\ord_\pp(x_{\ell_i}-x_1) = c + i$,
where the integer $c=\ord_\pp(x_{M+1}-x_1)$
is independent of $i$.
Therefore, for $i,j,k \in \Z_{\ge 1}$,
we have
\[
    i+j=k
    \quad\iff\quad
    \ord_\pp(x_{\ell_i}-x_1) + \ord_\pp(x_{\ell_j}-x_1)
    = \ord_\pp(x_{\ell_k}-x_1) + c.
\]
It follows that the graph of $+$ corresponds under $\phi$
to a subset of $A^3$ that is Diophantine over $\OO_{K,\calS}$.

Thus $A$ is a Diophantine model of $\calZ$ over $\OO_{K,\calS}$.
\end{proof}

As already remarked,
Corollary~\ref{C:positive existential model} and
Proposition~\ref{P:model of calZ}
together imply part~(2) of Theorem~\ref{thm:main2}.

\section*{Acknowledgements}

We thank Ernie Croot for a discussion
regarding Lemma~\ref{L:dirichlet-vinogradov}.



\providecommand{\bysame}{\leavevmode\hbox to3em{\hrulefill}\thinspace}

\end{document}